# Generalization of the Chekanov theorem: Diameters of immersed manifolds and wave fronts

Pushkar' P. E.

Moscow State University

e-mail – petya@pushkar.mccme.rssi.ru

## Introduction

The Chekanov theorem [6] generalizes the classic Lyusternik-Shnirel'man and Morse theorems concerning critical points of a smooth function on a closed manifold. A Legendrian submanifold $\Lambda$ of space of 1-jets of the functions on a manifold $M$ defines a multi-valued function whose graph is the projection of $\Lambda$ in $J^0M = M \times \mathbf{R}$. The Chekanov theorem asserts that if $\Lambda$ is homotopic to the 1-jet of a smooth function in the class of embedded Legendrian manifolds, then such a graph of a multi-valued function must have a lot of points (their number is determined by the topology of $M$) at which the tangent plane to the graph is parallel to $M \times 0$.

In the present paper a similar estimate is proved for a wider class of Legendrian manifolds. We consider Legendrian manifolds homotopic (in the class of embedded Legendrian manifolds) to Legendrian manifolds specified by generating families. Another generalization of the Chekanov theorem can be found in [5], [8], [9], [13].

As in [9], the proposed generalization of the Chekanov theorem is applied to investigating the geometry of wave fronts. The immersed manifold $M$ in $\mathbf{R}^n$ determines a Legendrian submanifold in the space $ST^*\mathbf{R}^n$ of the spherization of the cotangent bundle of $\mathbf{R}^n$ which can be defined by a generating family. An implicit form of this generating family makes it possible to obtain lower-bound estimates for the number of diameters (the diameter is a segment connecting different points of $M$ and orthogonal to $M$ at its endpoints) of the immersed manifolds $M$. These estimates comparing with the results of [12]. It is proved that the number of diameters (counted with multiplicities) is minorated by $\frac{1}{2}(B^2 + (\dim M - 1)B)$, where $B = \sum \dim H_*(M, \mathbf{Z}_2)$. Generalization of the Chekanov theorem makes it possible to extend these estimates to a wider class of hypersurfaces in $\mathbf{R}^n$ — wave fronts obtained from $M$ by means of a deformation without (dangerous) self-contact.

In §1 the generalization of the Chekanov theorem is proved. In §2 the number of points of self-intersection of the projection in $T^*M$ of a Legendrian submanifold (of the class considered) of the space $J^1M$ is estimated. In §3 the number of diameters of an immersed submanifold is estimated. In §4 the estimates of §3 are extended to wave fronts. In §5 the accuracy of the estimates of §3 and §4 is discussed.

# 1 Generalization of the Chekanov theorem. Critical points of quasifunctions

## 1.1 Formulation of the basic theorem

Let $p : E \to M$ be a bundle. The generic function $h$ on $E$ determines a Legendrian manifold in the space $J^1M$ . The manifold of critical points over the fiber can naturally be mapped into $J^1M$ (the differential $h$ "along the base" and the value of the function $h$ are related to the critical point). The pairs obtained form the Legendrian manifold $\Lambda \subset J^1M$.

*Definition.* A *$E$-quasifunction* is a Legendrian manifold homotopic to $\Lambda$ in $J^1M$ in the class of Legendrian embeddings.

*Definition* : *Critical points* of the Legendrian manifold $\Lambda$ in $J^1M$ are such points of $\Lambda$ whose images belong to the zero section under the natural mapping $\rho_M : J^1M \to T^*M$. Critical points correspond to the points of the front (graph of the $E$-quasifunction) at which the tangent plane is parallel to $M \times 0$. We will call a critical point non-degenerate, if $\rho_M(\Lambda)$ is transversal to the zero section at the image of the critical point. We now formulate the basic theorem.

**Theorem 1.1** *Suppose $\Lambda$ is a multi-valued $E$-quasifunction on a closed manifold $M$ and the bundle fiber$E \to M$ is compact. Then the number of critical points of $\Lambda$ is at least*

$\sum_{i=0}^{\dim E} (b_i(E) + 2q_i(E))$, *where $b_i(E) = \dim H_i(E, \mathbf{R})$ and $q_i(M)$ is the minimum number of generators of the group* Tors $H_i(E, \mathbf{Z})$, *if the critical points of $\Lambda$ are not degenerate.*

$\sum b_i(E, F)$, *$F$ is the field if the critical points of $\Lambda$ are not degenerate.*

$\mathrm{cl}(E, A) + 1$. *Here,* $\mathrm{cl}(E, A)$ *is the cohomology length of the manifold $E$ with coefficients in a commutative ring $A$ ([6]).*

## 1.2 Necessary notation.

We denote the canonical mapping $T^*B \to B$ by $\pi_B$ and the projection $J^1B \to T^*B$ (forgetting the value) by $\rho_B$.

Suppose $p : E \to M$ is a smooth bundle with the fiber $W$, $E_0 \subset T^*E$ is a sub-bundle of the bundle $\pi_E : T^*E \to E$ formed by covectors whose restriction to the tangent space of the fiber is equal to zero, and $\pi_\rho$ is the projection of $E_0$ in $E$. We denote the natural mapping of the bundle $E_0$ over $T^*M$ ($E_0 \times \mathbf{R}$ over $J^1M$) by $p_0$ ($p_1$).

The described mappings can be naturally unified in the commutative diagram:

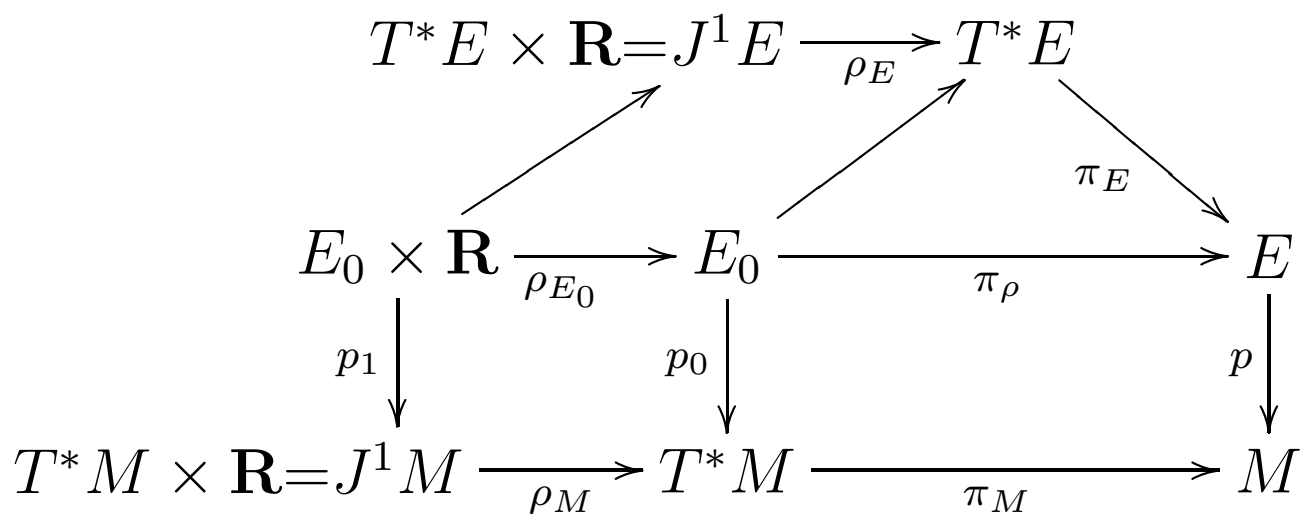

## 1.3 Definition of the generating family.

Suppose $p : E \to M$ is a bundle and $F : E \to \mathbf{R}$ is a smooth function so that the intersection $j^1F$ and $E_0 \times \mathbf{R}$ is transversal. Then $\Lambda = p_1\left(j^1F \cap (E_0 \times \mathbf{R})\right)$ is an (immersed) Legendrian manifold. When the transversality condition is satisfied, $F$ is called a *generating family* of the Legendrian manifold $\Lambda$. We will say that a generating family $\widetilde{F} : \widetilde{E} \to \mathbf{R}$ is a *stabilization* of a generating family $F : E \to \mathbf{R}$, if $\widetilde{E} = E \times \mathbf{R}^N$ and $\widetilde{F} = F + Q$, where $Q$ is non-degenerate quadratic form on $\mathbf{R}^N$. A generating family $\widetilde{F} : E \times \mathbf{R}^N \to \mathbf{R}$ will be called *quadratic at infinity* if it is the sum of a quadratic form on each fiber $\mathbf{R}^N$ and a function with bounded differential. We will sometimes denote the generating family $F$ by $F(x, q)$, where $q$ is a point of $M$ and $x$ is a point of the fiber $p$.

## 1.4 Proof of Theorem 1.1.

A function defined on the space of a vector bundle will be called *quadratic at infinity* if this function is the sum of a non-degenerate quadratic form along the fiber and a function with bounded differential. Then Theorem 1.1 directly follows from the following assertions:

**Theorem 1.2** *Suppose $M$ is a closed manifold and $p : E \to M$ bundle with a compact fiber. Then $E$-quasifunction $\Lambda$ on $M$ can be defined by the generating family $F : E \times \mathbf{R}^{N(\Lambda)} \to \mathbf{R}$ quadratic at infinity.*

**Theorem 1.3** *Suppose $V$ is a vector bundle over a closed manifold $E$ and $f : V \to \mathbf{R}$ is a function quadratic at infinity. Then the number of the critical points of the function $f$ is at least:*

$\sum\limits_{i=0}^{\dim E} (b_i(E) + 2q_i(M))$, *if all the critical point are Morse .*

$\sum b_i(E, F)$, *if all the critical point are the Morse and $F$ is a field of coefficients.*

$\mathrm{cl}(E, A) + 1$.

Theorem 1.3 follows from the well-known result obtained by Conley and Zehnder [7]. Theorem 1.2 can conveniently be proved in the following formulation:

**Theorem 1.4** *Suppose $\Lambda_0 \subset J^1M$ is a Legendrian manifold defined by a generating family $F : E \to \mathbf{R}$ and the fiber of the bundle $E$ is compact. Let $\{G^t\}_{t\in[0,1]}$ be a smooth family of contactomorphisms of the space $J^1M$, $G^0 = id$. Then $\Lambda_t = G^t(\Lambda_0)$ can be defined by the generating family $\widetilde{F_t} : E \times \mathbf{R}^N \to \mathbf{R}$ quadratic at infinity.*

*Proof:* Suppose $k_t$ is the contact Hamiltonian of the field $v_t = \left.\frac{\partial G^\tau}{\partial \tau}\right|_{t=\tau}$. Let us consider $K_t = p_1^* k_t : E_0 \times \mathbf{R} \to \mathbf{R}$ and extend $K_t$ to a function $\widetilde{K_t} : J^1E \to \mathbf{R}$. Then $E_0 \times \mathbf{R}$ is invariant with respect to the flow $\widetilde{G}^t$, given by the contact Hamiltonian $\widetilde{K_t}$,

$$p_1(\widetilde{G_t}(j^1F)) \cap (E_0 \times \mathbf{R}) = G^t(\Lambda_0) = \Lambda_t \qquad (*)$$

and $\widetilde{G}^t(j^1F)$ intersects $E_0 \times \mathbf{R}$ transversally. The function $\widetilde{K_t}$ can be improved so that it becomes finite, and $(*)$ and the transversality conditions are satisfied.

We can apply the Chekanov theorem ([6], Theorem 3.1) to the manifold $j^1F$ and the flow $\widetilde{G}^t$ determined by the contact Hamiltonian $\widetilde{K_t}$.

Consequently, $\widetilde{G}^t(j^1F)$ is defined by a generating family $F_t : E \times \mathbf{R}^{N(t)} \to \mathbf{R}$, quadratic at infinity. Accordingly, $G^t(\Lambda_0)$ is also defined by the generating family $F_t'$. The proof is completed. □

**Remark 1.5.** We can choose the number $N(t)$ to be independent on $t \in [0,1]$.

A Legendrian manifold $\Lambda_1 \subset J^1M$ will be called a *good $E$-quasifunction* on $M$, if $\Lambda_1$ is homotopic in the class of Legendrian embedding to a Legendrian manifold $\Lambda$ defined by a generating family $F : E \times \mathbf{R}^N \to \mathbf{R}$ quadratic at infinity.

**Remark 1.6.** Theorems 1.1 and 1.2 can be extended to good $E$-quasifunctions.

# 2 Lagrangian self-intersections.

Suppose $\Lambda$ is a Legendrian submanifold of $J^1M$ then $\rho_M(\Lambda)$ is an immersed Lagrangian submanifold of $T^*M$. If $\Lambda$ is defined by a generating family $F : E \to \mathbf{R}$ ($E$ is a bundle over $M$), then the number of self-intersection points $\rho_M(\Lambda)$ can be estimated in terms of the Betti numbers of $\Lambda$ and $E$.

**Theorem 2.1** *Suppose $\Lambda \subset J^1M$ is a $E$-quasifunction and $E \to M$ is a bundle with a compact fiber. Then the number of self-intersections points (counted with multiplicities) of the projection $\rho_M(\Lambda)$ is at least* $\frac{1}{2}\left|\sum b_i(\widetilde{E}, \mathbf{Z}_2) - \sum b_i(\Lambda, \mathbf{Z}_2)\right|$.

*Proof:* As follows from Theorem 1.2, $\Lambda$ can be defined by the generating family $F : E \times \mathbf{R}^N \to \mathbf{R}$, quadratic at infinity. We consider the bundle $\widetilde{E} \to M$ induced by the diagonal embedding of $M$ in $M \times M$ from the bundle $(E \times \mathbf{R}^N) \times (E \times \mathbf{R}^N) \to M \times M$. The fiber $\widetilde{E}$ is the Cartesian square of the bundle fiber $E$ multiplied by $\mathbf{R}^N$. On $\widetilde{E}$ we define the function $\widetilde{F}(x, y, q) = F(x, q) - F(y, q)$, where $x$ and $y$ are points of the bundle fiber $E \times \mathbf{R}^N \to M$, $q \in M$. The points of self-intersection of the projection are in one-to-one correspondence with the critical points of the function $\widetilde{F} : \widetilde{E} \to \mathbf{R}$ with a positive critical value. The points with zero critical value form a Bott manifold which is diffeomorphic to $\Lambda$ (see Lemma 2.3 below). In the general position the critical points $\widetilde{F}$ with nonzero critical value are Morse points. The number of the self-intersection points is half the number of the Morse points . By virtue Theorem 2.4 we obtain the required estimate. The proof is completed. $\square$

**Corollary 2.2** *Let $E = M \times W$ be a trivial bundle, and $\Lambda$ a multi-valued $E$-quasifunction. Then the number of points of self-intersection of the projection $\Lambda$ in $T^*M$ (in general position) is at least $\frac{1}{2}[\sum b_i(M) \left(\sum b_i(W)\right)^2 - \sum b_i(\Lambda)]$.*

**Lemma 2.3** *Suppose $F_1(x, q)$, $F_2(y, q)$ are generating families of a Legendrian manifold $\Lambda \subset J^1M$. Then critical points of the function $F_1(x, q) - F_2(y, q)$ with zero critical value form a Bott manifold diffeomorphic to $\Lambda$.*

*Proof:* The assertion that the manifold of critical points with zero critical value is diffeomorphic to the manifold $\Lambda$ is obvious. The assertion that this manifold is a Bott one is local. We will prove this assertion: Suppose $g_1(x, q)$ is a generating family of a germ of the Legendrian manifold $\Lambda$. The generating family $g_2 = g_1 + Q(z)$ will be called the stabilization of $g$, if $Q$ is a non-degenerate quadratic form. Two generating families $g_1(x, q)$, $g_2(x, q)$ of the germ of the Legendrian manifold $\Lambda$ are equivalently fibered if $g_1(x, q) = g_2(h(x, q), q)$, $(x, q) \to (h(x, q), q)$ is a diffeomorphism. Two stably fibered generating families will be called equivalent if these families become equivalent after stabilization. The assertion that the generating families $F_1$, $F_2$ are of the Bott type remains valid when $F_1$ and $F_2$ are replaced by stably fibered equivalent generating families. According to [3], if $g_1$ and $g_2$ are generating families of the same germ of a Legendrian manifold then $g_1$ is stably fibered equivalently to $g_2$. As $F_1$ and $F_2$ we choose the following generating family. Suppose that at the point considered the dimension of the kernel of the projection $\rho_m(\Lambda)$ in $M$ is $k$. Then there exist such canonical coordinates $q_1, \ldots, q_n$, $p_1, \ldots, p_n$ $T^*M$ that the kernel of the projection coincides with $\left\langle \frac{\partial}{\partial p_1}, \ldots, \frac{\partial}{\partial p_k} \right\rangle$ and the manifold $\rho_M(\Lambda)$ at the point $\rho_M(y)$ is tangent to $\left\langle \frac{\partial}{\partial p_1}, \ldots, \frac{\partial}{\partial p_k}, \frac{\partial}{\partial p_{k+1}}, \ldots, \frac{\partial}{\partial p_n} \right\rangle$. Then the germ of $\Lambda$ can be defined by the generating family $\sum_{i=1}^{k} x_i y_i + S(x_1, \ldots, x_k, q_{k+1}, \ldots, q_n)$, where $S = C + o(|x, q|^2)$. For this generating family we can readily show that this manifold is of the Bott type. □

We will need the following modification of Theorem 1.3.

**Theorem 2.4** *Suppose $V$ is the total space of a vector bundle over a closed manifold $M$, $g : V \to \mathbf{R}$ is a function quadratic at infinity, $L$ is a Bott manifold of the function $g$ lying at the level $g = C$, and the other critical points of $g$ are Morse points. Then the number of Morse critical points of $g$ is at least $|\sum b_i(M, \mathbf{Z}_2) - \sum b_i(L, \mathbf{Z}_2)|$.*

*Proof:* Without loss of generality, we can assume that all the critical values of $g$ are different. Let us consider the function

$$\alpha(t) = \dim H_*(g \leq C + a + t, g \leq C - t, \mathbf{Z}_2)$$

for a fairly small $a$ at $t > 0$.

From the Morse theory ([10], [4]) it follows that $\alpha(t)$ is changed by $\pm 1$ only when $C + a + t$ or $C - t$ pass through critical values. Consequently, the number of the critical points is not less than $|\alpha(t_1) - \alpha(t_0)|$, where $t_1$ is fairly large and $t_0$ is a fairly small positive number. The bundle $V$ can be decomposed into the sum of the bundles $V_+$ and $V_-$, and the function $g$ can be represented as the sum of a function on $V_+$ ($V_-$), namely, a fiber-wise positively (negatively) defined quadratic form, and a function with the bounded differential. From Morse theory it follows that $\alpha(t_1)$ is equal to the dimension of the space of $\mathbf{Z}_2$-homologies of the Thom space of the bundle $V_-$. From the Thom isomorphism we obtain that $\alpha(t_1) = \sum b_i(M, \mathbf{Z}_2)$. Analogous considerations for small $t$ show that $\alpha(t_0) = \sum b_i(L, \mathbf{Z}_2)$. The theorem is proved.

# 3 Minoration of the number of diameters of manifolds immersed in $\mathbf{R}^n$

## 3.1 Introduction

Suppose $f$ is an immersion of a manifold $M^n$ in $\mathbf{R}^{n+k}$. A segment connecting two different points $f(x)$ and $f(y)$ of the immersion and perpendicular to the tangent planes at these points will be called the *diameter* of an immersed manifold $f(M^n)$. We formulate the principal result:

**Theorem 3.1** *Suppose $M^n$ is a closed manifold of dimension $n$ and $B = \sum dim H_*(M, Z_2)$. Then for a generic immersion the number of diameters of $M^n$ in $\mathbf{R}^{n+k}$ is at least $\frac{1}{2}(B^2 + (n-1)B)$.*

The problem of diameters (double normals) of a single submanifold of a Riemannian manifold was considered in [], []. The case of non-general position for a sphere was known to Lyusternik, Shnirel'man and Morse. In the general position in Euclidean space the case of an immersed submanifold was considered by Takens and White in [12]. In that paper there was obtained an explicit estimate (for embeddings of $M^k$ in $\mathbf{R}^n$) which gives the number of the diameters in terms of Betti numbers $M$. We can show that this estimate does not exceed $\frac{1}{2}\left((\sum b_i(M))^2 - \sum b_i(M)\right) + 2n$. In all particular examples our estimate is not worse than the estimate [12]. (The case of a torus $T^n$ in [12] cannot be derived from other results of [12]. The author wishes to thank A. N. Mikoyan for his assistance in establishing this fact). We will give some corollaries.

**Corollary 3.2** *For generic immersion of a manifold $M$ in $\mathbf{R}^n$ the number $D$ of diameters $D$ is at least*
$M = S^n$, $D \geq n + 1$.
$M = T^n$, $D \geq 2^{2n-1} + (n-1)2^{n-1}$.
$M = S_g^2$ *(oriented surface of the kind $g$)*, $D \geq 2g^2 + 5g + 3$.
$M = \mathbf{R}P^n$, $D \geq n^2 + n$.

*Proof:* substitution to the formula of Theorem 3.1. □

## 3.2 Proof of Theorem 3.1.

**Remark.** The condition of generisity in Theorem 3.1 and below means that all the critical points of functions considered are Morse and self-intersections are transversal. We can show that this is actually the condition of the general position (like [12]).

Before to prove Theorem 3.1 we formulate certain convenient assertions.

We will consider function $F : S^{n+k-1} \times M \times M \to \mathbf{R}$, $F(x, \xi_1, \xi_2) = \langle x, f(\xi_1) - f(\xi_2)\rangle$. Here, $S^{n+k-1}$ is the sphere of radius 1 in $\mathbf{R}^{n+k}$ with its center at the origin , and $\langle ., .\rangle$ is a scalar product in $\mathbf{R}^{n+k}$.

**Lemma 3.3** *Suppose $f$ is a generic immersion and the point $(x', \xi_1', \xi_2')$ is critical for $F = \langle x, f(\xi_1) - f(\xi_2)\rangle$. If $f(\xi_1') \neq f(\xi_2')$ the segment $[f(\xi_1'), f(\xi_2')]$ is a diameter and $F(x', \xi_1', \xi_2') \neq 0$.*

*Proof:* Differentiating $F$ with respect to $\xi_i (i = 1, 2)$ at the point $(x', \xi_1', \xi_2')$, we obtain: $x' \perp f_*(T_{\xi_i'} M^n)$. Consequently, $f(\xi_1') \neq f(\xi_2')$ since the self-intersection is transversal. On the other hand, differentiating with respect to $x$, we obtain that $x'$ is proportional to $f(\xi_1') - f(\xi_2')$. Consequently, $f(\xi_1'), f(\xi_2')]$ is a diameter. Since $x'$ is proportional to $f(\xi_1') - f(\xi_2')$ then $F(x', \xi_1', \xi_2') \neq 0$. $\Box$

Thus, we construct a mapping of the set of critical points of $F$ into diameters. We can readily verify that this is the mapping surjective and that each diameter corresponds to exactly four critical points of $F$ (in the case of general position). Namely, the following assertion holds:

**Lemma 3.4** *Suppose $[f(\xi_1'), f(\xi_2')]$ is a diameter of $f(M^n)$ then*

$$\left(\pm\frac{f(\xi_1') - f(\xi_2')}{\|\xi_1 - \xi_2\|}, \xi_1, \xi_2\right), \left(\pm\frac{f(\xi_1') - f(\xi_2')}{\|\xi_1 - \xi_2\|}, \xi_2, \xi_1\right)$$

*are critical points of the function $F = \langle x, f(\xi_1) - f(\xi_2)\rangle$.*

*Proof of Theorem 3.1.* The function $F(x, \xi_1, \xi_2)$ is invariant with respect to the action of the involution $S^{n+k-1} \times M^n \times M^n$, $(x, \xi_1, \xi_2) \to (-x, \xi_2, \xi_1)$. The involution acts without fixed points, consequently, the quotient set of this involution is a smooth manifold (we denote this manifold as $S^{n+k-1} \times M^n \times M^n / \mathbf{Z}_2$). The function $F$ determines on the quotient a smooth function $\widetilde{F}$. From 3.3 and 3.4, it follows that the number of the diameters is not less than half the number of critical points with nonzero critical value. The critical points with nonzero critical value form a Bott manifold diffeomorphic to $P\nu(M^n)$, which is the projectivization of the normal bundle (see Lemma 3.5). Thus, from Morse theory it follows that $2D + \sum b_i(P\nu(M^n), \mathbf{Z}_2) \geq \sum b_i \left( S^{n+k-1} \times M^n \times M^n \big/_{\mathbf{Z}_2}, \mathbf{Z}_2 \right)$, where $D$ is the number of the diameters. Set $B = \sum b_i(M^n, \mathbf{Z}_2)$. From Lemma 3.6 it follows that

$$\sum b_i \left( S^{n+k-1} \times M^n \times M^n \big/_{\mathbf{Z}_2}, \mathbf{Z}_2 \right) = B^2 + (n+k-1)B$$

Lemma 3.7 implies that $\sum b_i(P\nu(M^n), \mathbf{Z}_2) \leq kB$. Thus, $D \geq \frac{1}{2}(B^2 + (k-1)B)$, which proves the theorem. □

**Lemma 3.5** *For an generic immersion the critical points of $\widetilde{F}$ with nonzero critical value form a Bott manifold diffeomorphic to the projectivization of the normal bundle.*

**Remark.** Here, generisity condition is the transversality of self-intersection.
*Proof:* From Lemma 3.3 it follows that the critical points with nonzero critical value of the function $F$ are contained in $S^{n+k-1} \times \Delta \subset S^{n+k-1} \times M^n \times M^n$ ($\Delta$ is the diagonal in $M^n \times M^n$). These points form a Bott manifold diffeomorphic to the spherization of the normal bundle. A point $(x, \xi_1, \xi_2)$ is critical if and only if $x \perp f_*(T_\xi(M^n))$. The fact that this manifold is Bott can be verified in the local coordinates. □

**Lemma 3.6** *Suppose $M$ is a closed manifold and $l$ is an involution on $S^N \times M \times M$, $l(x, y, z) = (-x, z, y)$. Then the sum of the $Z_2$-Betti number factor space of $S^N \times M \times M$ with respect to the action $l$ is equal to $B^2 + NB$.*

*Proof:* Suppose $(a_i)$ is a cell division of $M$. $\sigma^0, \widehat{\sigma}^0, \sigma^1, \widehat{\sigma}^1, \ldots, \sigma^N, \widehat{\sigma}^N$ is the standard cell division of $S^N$ invariant with respect to the antipodal involution. $\sigma^i, \widehat{\sigma}^i$ are the cells of the dimension $i$ and $l(\sigma^i) = \widehat{\sigma}^i$. Then $S^N \times M \times M$ is decomposed to the cells $\sigma^k \times \alpha^i \times \alpha^j$ and $\widehat{\sigma}^k \times \alpha^i \times \alpha^j$. The involution $l$ transfers the cell $\sigma^k \times \alpha^i \times \alpha^j$ to the cell $\widehat{\sigma}^k \times \alpha^i \times \alpha^j$. This cell division $S^N \times M \times M$ induces a cell division in $S^N \times M \times M \big/_l$ since it is invariant with respect to the involution.

Suppose $\partial$ is a boundary operator in $C_*(M)$, $\partial_1$ is a boundary operator in $C_*(M \times M)$, and $\partial_2$ is a boundary operator in $C_*(S^N \times M \times M\big/_l)$ ($\mathbf{Z}_2$ is the region of the coefficients). We identify (formally) the cells $S^N \times M \times M\big/_l$ with the cells $\sigma^k \times \alpha^i \times \alpha^j$. Then $\partial_2(\sigma^k \times \alpha^i \times \alpha^j) = \sigma^k(\partial_1(\alpha^i \times \alpha^j)) + \sigma^{k-1}(\alpha^i \times \alpha^j + \alpha^j \times \alpha^i)(\sigma^{-1} = 0)$.

Let us consider the case $\partial = 0$. In this case the complex of the cell chains $C_*(S^N \times M \times M\big/_l)$ is graded by the dimension of the cell in the sphere $S^N$, since from $\partial = 0$ there follows $\partial_1 = 0$. We denote homologies with respect to this graduation by $H_k$. We now calculate the dimension of $H_k$. In this graduation the complex $C_*(S^N \times M \times M\big/_l, Z_2)$ is as follows: for any $k, 0 \leq k \leq N$ $C_*(S^N \times M \times M\big/_l, Z_2)$ is isomorphic to $C_*(M \times M, Z_2)$ , and the boundary operator is the mapping of symmetrization $s$, $s(a_i \times a_j) = a_i \times a_j + a_j \times a_i$. $H_N = ker(s)$, consequently, $dimH_N = \frac{1}{2}(B^2 + B)$. Dimension $H_k$ for $0 < k < N$ is equal $B$, $dimH_0 = \frac{1}{2}(B^2 + B)$.

We now consider the case $\partial \neq 0$. We reduce it to the case $\partial = 0$. Using the given cell division $M$, we construct a cell space $M'$, so that the complex $C_*(S^N \times M \times M\big/_l, Z_2)$ and $C_*(S^N \times M' \times M'\big/_l, Z_2)$ are isomorphic and $dimH_*(S^N \times M' \times M'\big/_l, Z_2) = B^2 + NB$ . From this fact there follows the validity of the lemma in the case $\partial \neq 0$.

Any complex with the coefficients in the field can be decomposed into the sum of a (trivial) homology complex and two-term exact complexes.

We now consider the space $M''$ which is the union of spheres in which the number of the spheres of the dimension $k$ is equal to $dimH_k$. $M'$ is the wedge of $M''$ and the disks corresponding to short exact complexes. Then the spaces $(S^N \times M' \times M')\big/_l$ and $(S^N \times M'' \times M'')\big/_l$ are homotopically equivalent since $(S^N \times M'' \times M'')\big/_l$ is the strong retract of $(S^N \times M' \times M')\big/_l$, consequently, their homologies are identical. The case $(S^N \times M'' \times M'')\big/_l$ is considered in the case $\partial = 0$. $\Box$

**Lemma 3.7** *Suppose $E \to B$ is a bundle with the fiber $RP^k$ ($S^k$). Then $\sum b_i(E, \mathbf{Z}_2) \leq (k+1) \sum b_i(b, \mathbf{Z}_2)$ ( $\sum b_i(E, \mathbf{Z}_2) \leq 2 \sum b_i(B, \mathbf{Z}_2)$).*

*Proof:* The bundle $E \to B$ with the fiber $RP^k$ ($S^k$) is homotopically simple since the fiber homologies in any dimension are not greater than one-dimensional. Consequently, the spectral sequence (with the coefficients from $\mathbf{Z}_2$ calculates the homologies of bundle space $E$. The dimension of the term $E_2$ of this spectral sequence is equal to $\sum b_i(E,\mathbf{Z}_2) \leq (k+1)\sum b_i(b,\mathbf{Z}_2)$ ( $\sum b_i(E,\mathbf{Z}_2) \leq 2\sum b_i(B,\mathbf{Z}_2)$). The dimension of the homologies of $E$ is not greater than the dimension of the term $E_2$. The proof is completed.

# 4 Diameters of wave fronts

The aim of this section is to extend the domain of application of the estimates of §3.

We recall some standard facts of contact geometry.

Suppose $B$ is a smooth manifold. A (co-oriented) hyperplane in the tangent space at a given point is called a (co-oriented) contact element of the manifold $B$ applied at this point. All the (co-oriented) contact elements in $B$ form the space $PT^*B$ ($ST^*B$) of a projectivized cotangent bundle (spherization of the cotangent bundle). The space $PT^*B$ ($ST^*B$) has a contact structure defined canonically (see [2]).

An immersed manofold $X$ with transversal self-intersections in $B$ determines a Legendrian submanifold $P(X) \subset PT^*B$ ($S(X) \subset ST^*B$) representing a set of (co-oriented) contact elements tangent to $X$ (see [2]).

A (co-oriented) wave front in $B$ is the projection of a Legendrian submanifold $PT^*B$ ($ST^*B$) in $B$. For a generic Legendrian submanifold in $PT^*B$ ($ST^*B$) its (co-oriented) wave front in $B$ is the is a singular stratified (co-oriented) hypersurface which at any point has a (co-oriented) tangent plane. A generic Legendrian submanifold in $PT^*B$ ($ST^*B$) can be uniquely restored from its wave front in $B$. We will identify a Legendrian submanifold $PT^*B$ ($ST^*B$) with its wave front.

A (co-oriented) wave front $W_1$ in $B$ will be called a deformation without (dangerous) self-tangency of the front $W_0$ in $B$ if the Legendrian submanifold $PT^*B$ ($ST^*B$) corresponding to $W_1$ is a deformation of the Legendrian submanifold corresponding to $W_0$ in the class of embedded Legendrian submanifolds in $PT^*B$ ($ST^*B$).

A diameter of the wave front in $\mathbf{R}^{n+1}$ is called a segment connecting two different points of the wave front and perpendicular to it at its endpoints.

**Theorem 4.1** *Suppose a wave front L is a closed hypersurface with transversal self-intersections in* $\mathbf{R}^{n+1}$ *and* $L_1$ *is a deformation of L without self-tangencies. Then the wave front* $L_1$ *has at least* $\frac{1}{2}\left(\left(\sum b_i(L,\mathbf{Z}_2)\right)^2 + (n-1)\sum b_i(L,\mathbf{Z}_2)\right)$ *diameters counting with their multiplicity.*

Before proving Theorem 3.1 we formulate the following assertion.

**Proposition 4.2** *([1]) The space $ST^*\mathbf{R}^{n+1}$ is contactomorphic to $J^1S^n$.*

*Construction:* We identify $S^n$ with the standard unit sphere $||x|| = 1$ in $\mathbf{R}^{n+1}$ and the cotangent vector to $S^n$ with a vector perpendicular to $x$ (using the metric). Then to a point $(u,p,q)$ from $J^1S^n = \mathbf{R}\times T^*S^n$ we relate place a (co-oriented) contact element parallel to $x$ at the point $uq+p$. We obtain a mapping from $J^1S^n$ onto $ST^*\mathbf{R}^{n+1}$. We can verify that this is a contactomorphism. In what follows, we will identify $J^1S^n$ with $ST^*\mathbf{R}^{n+1}$ in the when need.

**Lemma 4.3** *Suppose $X$ is an immersed submanifold with transversal self-intersections in $\mathbf{R}^{n+1}$. Then $S(X)$ can be defined by the generating family $F : S^n\times X\to\mathbf{R}$, $F(q,x) = \langle q,x\rangle$.*

*Proof:* This is an assertion of the support function theory. In the Lagrangian case it can be found in [3].

Lemma 4.3 and Proposition 4.2 give many natural examples of Legendrian submanifolds in $J^1S^n$.

*Proof of Theorem 4.1* We will consider $L$ as a wave front of a Legendrian submanifold $S(L)\subset ST^*\mathbf{R}^{n+1}$. Then $L_1$ is a wave front of a $S^n\times L$-quasifunction symmetric with respect to the change in the co-orientation of a contact element. According to Theorem 1.2, this $S^n\times L$-quasifunction can be defined by the generating family $F : S^n\times L\times\mathbf{R}^N\to\mathbf{R}$ which is quadratic at infinity. In the general case the number of the diameters is less a quoter of the number of critical points of the function

$$F_1 : S^n\times L\times\mathbf{R}^N\times L\times\mathbf{R}^N\to\mathbf{R}$$

$$F_1(q,\xi_1,z_1,\xi_2,z_2) = F(q,\xi_1,z_1) + F(-q,\xi_2,z_2)$$
$$(q\in S^n, \xi_i\in L, z_i\in\mathbf{R}^N, i=1,2)$$

with nonzero critical value.

According to lemma 2.3, the manifold of the critical points with zero critical value is a Bott manifold since $-F(-x, \xi, z) = G(x, \xi, z)$ defines the same Legendrian submanifold as $F(x, \xi, z)$.

The function $F_1$ determines the function $\widetilde{F}_1$ on the factor $S^n \times L \times \mathbf{R}^N \times L \times \mathbf{R}^N$ according to the action of the involution $l_1 : l_1(q, \xi_1, z_1, \xi_2, z_2) = (-q, \xi_2, z_2, \xi_1, z_1)$. Critical points of $\widetilde{F}_1$ form a Bott manifold diffeomorphic to $L$. $S^n \times L \times \mathbf{R}^N \times L \times \mathbf{R}^N\big/_{l_1}$ is the vector bundle over $S^n \times L \times L\big/_{l}$ and the function $\widetilde{F}_1$ is quadratic at infinity.

Twice the number $2D$ of the diameters is equal to the number of the critical points. From Theorem 2.4 it follows that in the generic case

$$D \geq \frac{1}{2}\left(\sum b_i\left(S^n \times L \times L\big/_{l}, \mathbf{Z}_2\right) - \sum b_i(L, \mathbf{Z}_2)\right)$$

Using Lemma 3.6, we obtain the estimate required.

## 4.1 Passing and counterpassing diameters of co-oriented wave front.

Passing (counterpassing) diameter of a co-oriented wave front is a diameter of this wave front whose co-orientations at the endpoints have the same (opposite) directions.

**Theorem 4.4** *Suppose $L$ is a closed immersed submanifold in $\mathbf{R}^{n+1}$ with transversal self-intersections. Let $\Lambda_1$ be a deformation of $S(L)$ in the class of Legendrian submanifolds in $ST^*\mathbf{R}^{n+1}$. Then the co-oriented wave front $\Lambda_1$ has at least*

$(\sum b_i(L, \mathbf{Z}_2))^2 - \sum b_i(L, \mathbf{Z}_2)$ *passing diameters and*

$(\sum b_i(L, \mathbf{Z}_2))^2 + n \sum b_i(L, \mathbf{Z}_2)$ *counterpassing diameters.*

*The number of diameters of the wave front $\Lambda_1$ is at least*
$\frac{1}{2}\left((\sum b_i(\Lambda_1, \mathbf{Z}_2))^2 + n \sum b_i(\Lambda_1, \mathbf{Z}_2)\right)$.

*In all this case the diameters are counted with their multiplicities.*

*Proof:* a) $\Lambda_1$ is diffeomorphic to $S(L)$ and, according to Theorem 1.2 $\Lambda_1$ is a $S^n \times L$-quasifunction. Passing diameters of $\Lambda_1$ are in the one-to-one (in the general position) correspondence with the self-intersection points of $\rho_{S^n}(\Lambda_1)$.

Consequently, in the general position we can use the estimate of Corollary 2.2.

b) $\Lambda_1$ is a $S^n \times L$-quasifunction. Consequently, in accordance with Theorem 1.2, $\Lambda_1$ can be defined by the generating family $F : S^n \times L \times \mathbf{R}^N \to \mathbf{R}$ quadratic at infinity. Counterpassing diameters are the critical points of the function $F_1$ defined in the proof of Theorem 4.1. In the general position all the critical points are non-degenerate and their number can be estimated using Lemma 3.6 and Morse theory.

c) This assertion can be proved if we sum estimates a) and b) and take into account that $\sum b_i(\Lambda_1, \mathbf{Z}_2) = \sum b_i(S(L), \mathbf{Z}_2) \leq 2 \sum b_i(L, \mathbf{Z}_2)$ (according to Lemma 3.7).

**Remark.** The assertion c) is not a corollary of Theorem 4.1, since in Theorem 4.1 we forbid all self-contacts, where as in c) only dangerous ones.

**Remark.** For the case of $L$ being a point Theorem 4.4 b) was proved in [9].

# 5 Sharpness of the estimates of the number of diameters of wave fronts

The aim of this section is to discuss the estimates of §4 and §3.

The estimate of §4 is correct for the $n$-dimensional sphere $S^n$ realized by an ellipsoid with different semiaxes in $(n+1)$-dimensional space.

Let us consider the Euclidean space $\mathbf{R}^{n+k+1}$ and the standard sphere $S^n$ of unit radius in $\mathbf{R}^{n+1} \subset \mathbf{R}^{n+k+1}$. Suppose $L$ is the set of those points of $\mathbf{R}^{n+k+1}$, the distance from which to $S^n$ is equal to $\frac{1}{2}$. Suppose that $L$ is arbitrarily oriented.

**Theorem 5.1** *The manifold $L$ is diffeomorphic to $S^n \times S^k$ and we can move $L$ and reduce it to the general position so that the number of the diameters will be equal to $2(n+k)+6$. Then the numbers of the passing and counterpassing diameters are equal to $2$ and $2(n+k)+4$.*

In fact, $L$ is diffeomorphic to $S^n \times S^k$ since the normal bundle to $S^n$ can be made trivial. We consider the function $f = \|\xi_1 - \xi_2\|^2 \ \ (L \times L \setminus \Delta)/_i$ ($i$ is the transposition of the factors). We can readily verify that the function $f$ on $(L \times L \setminus \Delta)/_i$ has four Bott manifolds. One of them is diffeomorphic to $RP^k \times S^n$, two are diffeomorphic to $\mathbf{R}P^n$, and one is diffeomorphic to $S^n$.

We can move $L$ so that in the vicinity of these Bott manifolds $f$ would have exactly $\sum b_i$ of Morse critical points (this fact we do not prove here). In this case, $\mathbf{R}P^k \times S^n$ and two copies of $\mathbf{R}P^n$ correspond to the counterpassing diameters, and $S^n$ to the passing diameters.

**Acknowledgments.** I am very grateful to V.I. Arnol'd, Yu.V. Chekanov and A.G. Khovanskii for fruitful discussions. The work was partially supported by the grant INTAS-4373 and RBRF (grant 96-01-01104) . I am very thankful to The Fields Institute where the work was completed for its hospitality. The present note is an announcement of the main results of [11].